\documentclass{elsart}
\usepackage{amsmath,amscd,amssymb}

\allowdisplaybreaks
\numberwithin{equation}{section}

\numberwithin{thm}{section}
\def\qedsymbol{~\qed}
\def\proofname{Proof}
\newenvironment{proof}{\par\noindent{\it\proofname}.~}{{\unskip\nobreak\hfill{\it\qedsymbol}}\par\vskip 9pt}
\newenvironment{proof*}[1]{\par\noindent{\it#1}.~}{{\unskip\nobreak\hfill{\it\qedsymbol}}\par\vskip 9pt}
\newtheorem{ass}{Assumption}
\newcommand{\cat}[1]{\operatorname{cat}{(#1)}}
\newcommand{\Cat}[1]{\operatorname{Cat}{(#1)}}
\newcommand{\wcat}[1]{\operatorname{{\it w}cat}{(#1)}}
\newcommand{\sigmacat}[1]{\operatorname{{\it r}cat}{(#1)}}

\newcommand{\cuplen}{\operatorname{cup}}

\newcommand{\weight}{\operatorname{wgt}}

\begin{document}
\begin{frontmatter}

\title{Lusternik-Schnirelmann category of non-simply connected compact simple Lie groups
}
\author[Kyushu]{Norio Iwase\thanksref{X}}
\ead{iwase@math.kyushu-u.ac.jp}
\author[Okayama]{Mamoru Mimura}
\ead{mimura@math.okayama-u.ac.jp}
\author[Kinki]{Tetsu Nishimoto}
\ead{nishimoto@kinwu.ac.jp}
\address[Kyushu]{Faculty of Mathematics,
  Kyushu University,
  Ropponmatsu Fukuoka 810-8560, Japan}
\address[Okayama]{Department of Mathematics,
  Faculty of Science,
  Okayama University,
  3-1 Tsushima-naka, Okayama 700-8530, Japan}
\address[Kinki]{Department of Welfare Business,
  Kinki Welfare University, Fukusaki-cho,
  Hyogo 679-2217, Japan}

\thanks[X]{The first named author is supported by the Grant-in-Aid for
  Scientific Research \#14654016 from Japan Society for the Promotion
  of Science.}

\begin{abstract}
Let $F \hookrightarrow X \to B$ be a fibre bundle with structure group $G$, where $B$ is $(d{-}1)$-connected and of finite dimension, $d \geq 1$.
We prove that the strong L-S category of $X$ is less than or equal to $m + \begin{textstyle}\frac{\dim B}{d}\end{textstyle}$, if $F$ has a cone decomposition of length $m$ 
under a compatibility condition with the action of $G$ on $F$.
This gives a consistent prospect to determine the L-S category of non-simply connected Lie groups.
For example, we obtain $\cat{PU(n)} \leq 3(n{-}1)$ for all $n \geq 1$, which might be best possible, since we have $\cat{\mathrm{PU}(p^r)}=3(p^r{-}1)$ for any prime $p$ and $r \geq 1$.
Similarly, we obtain the L-S category of $\mathrm{SO}(n)$ for $n \leq 9$ and $\mathrm{PO}(8)$.
We remark that all the above Lie groups satisfy the Ganea conjecture on L-S category.
\end{abstract}

%
%

\end{frontmatter}

\section{Introduction}

The Lusternik-Schnirelmann category $\cat{X}$, L-S category for short, is the least integer $m$ such that there is a covering of $X$ by $(m{+}1)$ open subsets each of which is contractible in $X$.

Ganea \cite{ganea} introduced a stronger notion of L-S category, $\Cat{X}$, which is equal to the cone-length, that is, the least integer $m$ such that there is a set of cofibre sequences $\{A_{i} \to X_{i-1} \hookrightarrow X_{i}\}_{1 \leq i \leq m}$ with $X_{0} = \{\ast\}$ and 
$X_{m}$ homotopy equivalent to $X$.

The weak L-S category $\wcat{X}$ is the least integer $m$ such that the reduced diagonal map $\bar\Delta^{m+1} : X \rightarrow {\bigwedge}^{m+1}X$ is trivial where ${\bigwedge}^{m+1}X$ is the smash product.
The stabilised version of the invariant $\wcat{X}$ is given as the least integer $m$ such that the reduced diagonal map $\bar\Delta^{m+1} : X \rightarrow {\bigwedge}^{m+1}X$ is {\it stably} trivial.
Let us denote it by $\cuplen(X)$, the {\it cup-length} of $X$.

In 1971, Ganea \cite{Ganea:conjecture} posed 15 problems on L-S category and its related topics: Computation of L-S category for various manifolds is given as the first problem and the second problem is known as the Ganea conjecture on L-S category.
These problems especially the first two problems have attracted many authors such as James and Singhof \cite{James:ls-category}, \cite{Singhof:minimal-cover}, \cite{Singhof:cat-lie1}, \cite{Singhof:cat-lie2}, \cite{Singhof:cat-lie3}, \cite{JS:SO(5)}, G\'omez-Larra\~naga and Gonz\'alez-Acu\~na \cite{GG:ls-cat_3m}, Montejano \cite{Montejano:singhof}, Oprea and Rudyak \cite{Rudyak:ls-cat_mfds1}, \cite{Rudyak:ls-cat_mfds2}, \cite{OR:ls-cat_mfds} and the authors \cite{Iwase:counter-ls}, \cite{Iwase:counter-ls-m}, \cite{Iwase:ls-cat-bundle}, \cite{IM:cat-sp(3)}, \cite{IMN:spin(7)}.
In \cite{Iwase:counter-ls-m,Iwase:ls-cat-bundle}, the first author gave a counter example as a manifold to the Ganea conjecture on L-S category.

Especially for L-S category of compact connected simple Lie groups,
the followings have already been known:
\begin{align*}
  &
  \cat{\mathrm{Sp}(1)}=\cat{\mathrm{SU}(2)}=\cat{\mathrm{Spin}(3)}=1,
  \\
  & \cat{\mathrm{SU}(3)}=2, \hspace{50pt}
  \cat{\mathrm{SO}(3)}=3,
\end{align*}
since $\mathrm{Sp}(1)=\mathrm{SU}(2)=\mathrm{Spin}(3) = S^3$,
$\mathrm{SU}(3) = {\Sigma}{\mathbb C}P^2 \cup e^8$ and $\mathrm{SO}(3)
= {\mathbb R}P^3$.
Schweitzer \cite{sch} showed
\begin{displaymath}
  \cat{\mathrm{Sp}(2)}=3
\end{displaymath}
using functional cohomology operations.
Singhof \cite{Singhof:cat-lie1,Singhof:cat-lie3} showed 
\begin{align*}
  & \cat{\mathrm{SU}(n)}=n{-}1, \\
  & \cat{\mathrm{Sp}(n)} \geq n+1, \quad \textrm{if } n \geq 2.
\end{align*}
Also we know
\begin{displaymath}
  \cat{\mathrm{G}_2}=4
\end{displaymath}
by \cite{James:ls-category} (see \cite{IM:cat-sp(3)}).
James and Singhof \cite{JS:SO(5)} showed
\begin{displaymath}
  \cat{\mathrm{SO}(5)}=8.
\end{displaymath}
The first and second authors \cite{IM:cat-sp(3)} and
Fern\'andez-Su\'arez, G\'omez-Tato, Strom and Tanr\'e \cite{fgst}
proved
\begin{align*}
  & \cat{\mathrm{Sp}(3)}=5, \\
  & \cat{\mathrm{Sp}(n)} \geq n+2 \quad \textrm{if } n \geq 3,
\end{align*}
by showing the reduced diagonal $\bar{\Delta}^5$ is given by the Toda
bracket $\{\eta,\nu,\eta\}=\nu^2$.
The authors \cite{IMN:spin(7)} showed
\begin{displaymath}
  \cat{\mathrm{Spin}(7)}=5, \quad \cat{\mathrm{Spin}(8)}=6
\end{displaymath}
using explicit cone decompositions of $\mathrm{Spin}(7)$ and $\mathrm{SU}(4)$.
Then the Ganea conjecture on L-S category holds for all these Lie groups, since the L-S and the strong L-S categories are equal to the cup-length:
\begin{fact}
If $\cat{X}=\cuplen{X}$, then the Ganea conjecture on L-S category holds for $X$, i.e., $\cat{X{\times}S^n}=\cat{X}{+}1$ for all $n \geq 1$.
\end{fact}
In fact, we have $\cuplen(X \times S^n) = \cuplen(X){+}1$ in general.

For any multiplicative cohomology theory $h$, we define $\cuplen(X;h)$, the {\it cup-length} with respect to $h$, by the least integer $m$ such that $u_0{\cdot} \cdots {\cdot}u_m = 0$ for any $m{+}1$ elements $u_i \in \tilde{h}^*(X)$.
When $h$ is the ordinary cohomology theory with coefficient ring $R$, $\cuplen(X;h)$ is often denoted as $\cuplen(X;R)$.
\begin{thm}
For any CW-complex $X$ we have
$$
\cuplen(X) = \max \{ \cuplen(X;h) \mid \textrm{$h$ is any multiplicative cohomology theory} \}.
$$
\end{thm}
\begin{proof}
It is easy to see that $\cuplen(X) \geq \cuplen(X;h)$, and hence we have $\cuplen(X)$ $\geq$ $\max \{ \cuplen(X;h) \mid \text{$h$ is any multiplicative cohomology theory} \}$.
Thus we must show 
$$
\cuplen(X) \leq \max \{ \cuplen(X;h) \hspace{-1pt} \mid \hspace{-1pt} \text{$h$ is any multiplicative cohomology theory} \}.
$$
Let $m = \max \{ \cuplen(X;h) \mid \text{$h$ is any multiplicative cohomology theory} \}$ and $h_X$ be the multiplicative cohomology theory represented by the following wedge sum of iterated smash products of suspension spectrum ${\Sigma}^{\infty}X$:
\begin{equation*}\begin{textstyle}
  S^0 \vee {\Sigma}^{\infty}X \vee {\Sigma}^{\infty}{\bigwedge}^{2}X \vee \cdots \vee {\Sigma}^{\infty}{\bigwedge}^{i}X \vee \cdots.
\end{textstyle}\end{equation*}
Let $\iota \in \tilde{h}^*_X(X)$ be the element which is represented by the inclusion map into the second factor ${\Sigma}^{\infty}X$ of the above wedge sum.
Then by the definition of the cup-length, we have $\iota^{m{+}1} = 0$ which is represented by the reduced diagonal map $\bar\Delta^{m{+}1} : X \to \bigwedge^{m{+}1}X$ in the $(m{+}2)$-nd factor ${\Sigma}^{\infty}\bigwedge^{m{+}1}X$ of the above wedge sum.
Hence we have $\cuplen(X) \leq m$ the desired inequality.
Thus we obtain the result.
\end{proof}
Let $P^{m}(\Omega{X})$ be the $m$-th projective space, in the sense of Stasheff \cite{Stasheff:higher-associativity}, such that there is a homotopy equivalence $P^{\infty}(\Omega{X}) \simeq X$.
The following theorem is obtained by Ganea (see also \cite{Iwase:counter-ls} and Sakai \cite{Sakai:push_pull}).
\begin{thm}[Ganea \cite{ganea}]
$\cat{X} \leq m$ if and only if there is a map $\sigma : X \to P^{m}(\Omega{X})$ such that $e_m^{X}{\circ}\sigma \sim 1_{X}$, where $e_m^{X} : P^{m}(\Omega{X}) \hookrightarrow P^{\infty}(\Omega{X}) \simeq X$.
\end{thm}
Using this, Rudyak \cite{Rudyak:ls-cat_mfds2,Rudyak:ls-cat_mfds3} introduced a stable L-S category, $\sigmacat{X}$, which is the least integer $m$ such that there is a stable map $\sigma : X \to P^{m}(\Omega{X})$ satisfying $e_m^{X}{\circ}\sigma \sim 1_{X}$, another stabilised version of L-S category.

Rudyak \cite{Rudyak:ls-cat_mfds1} \cite{Rudyak:ls-cat_mfds2} and Strom \cite{strom} introduced the following invariant to calculate $\sigmacat{X}$: Let $\weight(X;h)$ be the least integer $m$ such that the homomorphism $(e_m^{X})^{\ast} : \tilde{h}^{\ast}(X) \to \tilde{h}^{\ast}(P^{m}(\Omega{X}))$ is injective for any cohomology theory $h$.
When $h$ is the ordinary cohomology theory with coefficient ring $R$, $\weight(X;h)$ is often denoted as $\weight(X;R)$.

Since a product of any $m{+}1$ elements of $\tilde{h}^{\ast}(P^{m}(\Omega{X}))$ is trivial, we have $\cuplen(X;h) \leq \weight(X;h)$ for any multiplicative cohomology theory $h$.
Hence we have $\cuplen(X) \leq \weight(X)$, where we denote $\weight(X) = \max \{ \weight(X;h) \mid \text{$h$ is any cohomology theory}\}$.
\begin{rem}
For any ring $R$, we know $\cuplen(\mathrm{Sp}(2);R) = \weight(\mathrm{Sp}(2);R) = 2 < 3 = \cat{\mathrm{Sp}(2)}$.
But an easy calculation of algebra structure of $KO^*(\mathrm{Sp}(2))$ yields $\cuplen(\mathrm{Sp}(2);KO) = \weight(\mathrm{Sp}(2);KO) = 3 = \cat{\mathrm{Sp}(2)}$.
\end{rem}

The following theorem is due to Rudyak \cite{Rudyak:ls-cat_mfds2,Rudyak:ls-cat_mfds3}, although we do not know the precise relation between $\wcat{X}$ and $\sigmacat{X}$.
\begin{thm}
For any CW complex $X$, we have 
\begin{equation*}
\sigmacat{X} = \weight{X}
\end{equation*}
and hence we have the following relations among categories:
\begin{equation*}
  \cuplen(X) \leq \wcat{X}, \sigmacat{X}
  \leq \cat{X} \leq \Cat{X}.
\end{equation*}
\end{thm}
Using this stabilised version of L-S category, we have the following theorem.
\begin{thm}[Rudyak \cite{Rudyak:ls-cat_mfds2,Rudyak:ls-cat_mfds3}]\label{thm:rudyak}%
If $\cat{X}=\sigmacat{X}$, then the Ganea conjecture on L-S category holds for $X$.
\end{thm}
In fact, we have $\sigmacat{X \times S^n} = \sigmacat{X}{+}1$ in general (\cite{Rudyak:ls-cat_mfds2,Rudyak:ls-cat_mfds3}).

\section{Main results}\label{sec:main}

From now on, we work in the category of connected CW-complexes and continuous maps.
We denote by $Z^{(k)}$ the $k$-skeleton of a CW complex $Z$.
\begin{thm}[James \cite{James:ls-category},Ganea \cite{ganea}]\label{thm:james-ganea}
Let $X$ be a $(d-1)$-connected space of finite dimension.
Then $\begin{textstyle}\cat{X} \leq \Cat{X} \leq [\frac{\dim{(X)}}{d}]\end{textstyle}$, where $[a]$ denotes the biggest integer $\leq a$.
\end{thm}

In this paper, we extend this for a total space of a fibre bundle, to
determine L-S categories of $\mathrm{SO}(n)$ for $n \leq 9$,
$\mathrm{PO}(8)$ and $\mathrm{PU}(p^r)$ (and the other quotient groups
of $\mathrm{SU}(p^r)$), which also gives an alternative proof of a
result due to James and Singhof \cite{JS:SO(5)} on $\mathrm{SO}(5)$.

We assume that $B$ is a $(d{-}1)$-connected finite dimensional CW complex ($d\geq1$), whose cells are concentrated in dimensions $0,1,\cdots,s \mod d$ for some $s$, ($0 {\leq} s {\leq} d{-}1$).
Let $F \hookrightarrow X \to B$ be a fibre bundle with structure group $G$, a compact Lie group.
Then we have the associated principal bundle $G \hookrightarrow E \overset{\pi}{\rightarrow} B$ with $G$-action $\psi : G {\times} F \to F$ on $F$ and hence $X = E{\times}_GF$.

Let $K_{i} \overset{\rho_{i}}\to F_{i{-}1} \hookrightarrow F_{i}$, (${1 {\leq} i {\leq} m}$) be $m$ cofibre sequences with $F_{0} = \{\ast\}$ and $F_{m}$ homotopy equivalent to $F$.
We consider the following compatibility condition of the above cone decomposition of $F$ and the action of $G$ on $F$.
\begin{ass}\label{ass:compatible-decomposition}
$\psi \vert_{G^{(d{\cdot}(i{+}1){+}s{-}1)}{\times}F_{j}} : G^{(d{\cdot}(i{+}1){+}s{-}1)}{\times}F_{j} \to F$ is compressible into $F_{i{+}j}$, $0 \leq i,j \leq i{+}j \leq m$.
\end{ass}
\begin{rem}\label{rem:cellular}
\begin{enumerate}
\item\label{rem:cellular1}
Let $F=G$ and $X=E$ be the total space of a principal bundle over a path-connected space $B$ and $d=1$.
Then any cone decomposition of $F$ such that $F_{i}=F^{(n_{i})}$ with $0 < n_{1} < n_{2} < \cdots < n_{m} = \dim(F)$
satisfies Assumption \ref{ass:compatible-decomposition} with $s=0$.
\item\label{rem:cellular2}
Let $F \hookrightarrow X \to B$ be a trivial bundle.
Then any cone decomposition of $F$ satisfies the compatibility Assumption \ref{ass:compatible-decomposition} with $s=d{-}1$.
\end{enumerate}
\end{rem}

Our main result is stated as follows:

\begin{thm}
  \label{thm:main}
Let $B$ be a $(d{-}1)$-connected finite dimensional CW complex ($d\geq1$), whose cells are concentrated in dimensions $0,1,\cdots,s \mod d$ for some $s$, $0 \leq s \leq d{-}1$.
Let $F \hookrightarrow X \to B$ be a fibre bundle with fibre $F$ whose structure group is a compact Lie group $G$.
If $F$ has a cone decomposition with the compatibility Assumption \ref{ass:compatible-decomposition} for $d$, then $\begin{textstyle}\Cat{X} \leq m + [\frac{\dim B}{d}]\end{textstyle}$.
\end{thm}

\begin{cor}\label{cor:main}
If $F$ has a cone decomposition with the compatibility Assumption \ref{ass:compatible-decomposition} for $s=d{-}1$ and also $m=\Cat{F}$, then $\begin{textstyle}\Cat{X} \leq \Cat{F} + [\frac{\dim B}{d}]\end{textstyle}$.
\end{cor}
\begin{rem}
Without Assumption \ref{ass:compatible-decomposition}, we only have 
$$\Cat{X}{+}1 \leq (\Cat{F}{+}1){\cdot}(\Cat{B}{+}1)$$
which is obtained immediately from the definition of $\mathrm{Cat}$ by Ganea \cite{ganea} and the corresponding results of Varadarajan \cite{Varadarajan:fib-cat} and Hardie \cite{Hardie:fib-cat} for $\mathrm{cat}$.
For example, the principal bundle $\mathrm{Sp}(1) \hookrightarrow \mathrm{Sp}(2) \to S^7$ does satisfy Assumption \ref{ass:compatible-decomposition} for $d \leq 3$, but not if $d \geq 4$, and we have $\Cat{\mathrm{Sp}(2)}\leq\Cat{\mathrm{Sp}(1)}{+}[\frac{7}{3}]=3>2=\Cat{\mathrm{Sp}(1)}{+}[\frac{7}{4}]$.
In fact by Schweitzer \cite{sch}, we know $\Cat{\mathrm{Sp}(2)}=3$.
\end{rem}
\begin{rem}
By Remark \ref{rem:cellular} (\ref{rem:cellular2}), Theorem \ref{thm:main} generalises Theorem \ref{thm:james-ganea}.
\end{rem}

By applying this, we first obtain the following general result:
\begin{thm}\label{thm:maing}
Let $C_m < SU(n)$ be a central (cyclic) subgroup of order $m$.
Then we have $\Cat{\mathrm{SU}(n)/C_{m}} \leq 3(n{-}1)$ for all $n \geq 1$.
\end{thm}
This might be best possible, because we also obtain the following result.
\begin{thm}
  \label{thm:main2}
  We have
  \begin{displaymath}
    \Cat{\mathrm{SU}(p^r)/C_{p^s}} = \cat{\mathrm{SU}(p^r)/C_{p^s}}
    = \sigmacat{\mathrm{SU}(p^r)/C_{p^s}} = 3(p^r{-}1)
  \end{displaymath}
  where $p$ is a prime and $1 \leq s \leq r$.
\end{thm}

Similarly we obtain the following result.

\begin{thm}
  \label{thm:main1}
  We have
  \begin{align*}
    & \Cat{\mathrm{SO}(6)} = \cat{\mathrm{SO}(6)} = \cuplen(\mathrm{SO}(6)) = 9, \\[2mm]
    & \Cat{\mathrm{SO}(7)} = \cat{\mathrm{SO}(7)} = \cuplen(\mathrm{SO}(7)) = 11, \\[2mm]
    & \Cat{\mathrm{SO}(8)} = \cat{\mathrm{SO}(8)} = \cuplen(\mathrm{SO}(8)) = 12, \\[2mm]
    & \Cat{\mathrm{SO}(9)} = \cat{\mathrm{SO}(9)} = \cuplen(\mathrm{SO}(9)) = 20, \\[2mm]
    & \Cat{\mathrm{PO}(8)} = \cat{\mathrm{PO}(8)} = \cuplen(\mathrm{PO}(8)) = 18.
  \end{align*}
\end{thm}
\begin{rem}
Theorem \ref{thm:main} also provides an alternative proof for a result of James-Singhof \cite{JS:SO(5)}, that is, $\Cat{\mathrm{SO}(5)} = \cat{\mathrm{SO}(5)}= \cuplen(\mathrm{SO}(5)) = 8$ (see Section \ref{sect:3}).
\end{rem}

We summarise all the known cases in the following table, where each number given in the right hand side of a connected, compact, simple Lie group indicates its L-S category.
\begin{equation*}
  \begin{array}{|c|c|c|c|c|c|c|c|c|c|c|}
    \hline
    \text{rank} & \multicolumn{2}{c}{1}
    & \multicolumn{2}{|c}{2}
    & \multicolumn{2}{|c}{3}
    & \multicolumn{2}{|c}{4}
    & \multicolumn{2}{|c|}{n~ (\geq 5)} \\
    \hline
    A_n     & \mathrm{SU}(2) & 1 & \mathrm{SU}(3) & 2 & \mathrm{SU}(4)   & 3  & \mathrm{SU}(5)   & 4
    & \mathrm{SU}(n{+}1) & n \\
            &       &   &       &   & \mathrm{SO}(6)   & 9  &         &    & \begin{array}{c}\vspace{-3pt}\vdots\end{array} & \\
            & \mathrm{PU}(2) & 3 & \mathrm{PU}(3) & 6 & \mathrm{PU}(4)   & 9  & \mathrm{PU}(5)   & 12 & \mathrm{PU}(n{+}1) & - \\
    \hline
    B_n     & \mathrm{Spin}(3) & 1 & \mathrm{Spin}(5) & 3 & \mathrm{Spin}(7) & 5  & \mathrm{Spin}(9) & - & \mathrm{Spin}(2n{+}1) & - \\
            & \mathrm{SO}(3) & 3 & \mathrm{SO}(5) & 8  & \mathrm{SO}(7)   & 11 & \mathrm{SO}(9)   & 20 & \mathrm{SO}(2n{+}1) & - \\
    \hline
    C_n     & \mathrm{Sp}(1) & 1 & \mathrm{Sp}(2) & 3 & \mathrm{Sp}(3)   & 5  & \mathrm{Sp}(4)   & - & \mathrm{Sp}(n) & - \\
            & \mathrm{PSp}(1) & 3 & \mathrm{PSp}(2) & 8 & \mathrm{PSp}(3)  & - & \mathrm{PSp}(4)  & - & \mathrm{PSp}(n) & - \\
    \hline
    D_n     &       &   &       &   & \mathrm{Spin}(6) & 3  & \mathrm{Spin}(8) & 6  & \mathrm{Spin}(2n) & - \\
            &       &   &       &   & \mathrm{SO}(6)   & 9  & \mathrm{SO}(8)   & 12 & \mathrm{SO}(2n) & - \\
            &       &   &       &   & \mathrm{PO}(6)   & 9  & \mathrm{PO}(8)   & 18 & \mathrm{PO}(2n) & - \\
            &       &   &       &   &         &    &         &    & \mathrm{Ss}(2n) & - \\
    \hline
    \text{\begin{tabular}{c} Except.\\[-2mm]types\end{tabular}} & &   & \mathrm{G}_2   & 4 &         &    & \mathrm{F}_4     & - & \mathrm{E}_6, \mathrm{E}_7, \mathrm{E}_8 & - \\
    \hline
  \end{array}
\end{equation*}
where "-" indicates the unknown case.
\begin{rem}
We recall that $A_1 = B_1 = C_1$, $B_2 = C_2$ and $A_3 = D_3$, and that the semi-spinor group $\mathrm{Ss}(2n)$ is defined only for $n$ even.
\end{rem}

Taking into account the above table, we get the following by Theorem \ref{thm:rudyak}:
\begin{cor}
The Ganea conjecture on L-S category holds for every connected, compact, simple Lie group $G$ when L-S category is known as above.
\end{cor}

The paper is organised as follows; In Section \ref{sect:2} we prove Theorem \ref{thm:main}.
In Section \ref{sect:3} we determine $\cat{\mathrm{SO}(n)}$ for $n = 5, 6, 7, 8, 9$ and $\cat{\mathrm{PO}(8)}$.
In Section \ref{sect:4} we prove Theorem \ref{thm:maing} and determine $\cat{\mathrm{SU}(p^r)/C_{p^s}}$.

\section{Proof of Theorem \ref{thm:main}}\label{sect:2}

Let $B_i$ be the $(d{\cdot}i{+}s)$-skeleton of $B$ and $n{=}[\frac{\dim B}{d}]$ the biggest integer not exceeding $\begin{textstyle}\frac{\dim B}{d}\end{textstyle}$.
Then by Ganea \cite{ganea}, Theorem \ref{thm:james-ganea} implies that there are $n$ cofibre sequences $A_{i} \overset{\lambda_{i}}\to B_{i{-}1} \hookrightarrow B_{i}$, $1 {\leq} i {\leq} n$ with $B_{0} = \{\ast\}$, $B_{n} = B$.
Note that $A_{i}$ is $(d{\cdot}i{-}2)$-connected and of dimension $(d{\cdot}i{+}s{-}1)$.
Hence we obtain
\begin{align*}&\begin{textstyle}
B_{i} = B_{i{-}1} \cup_{\lambda_{i}} C(A_{i}),\quad 
\lambda_i : A_{i} \to B_{i{-}1}\quad 
\end{textstyle}\\[2mm]&\qquad\quad\begin{textstyle}
A_{i} = A_{i}^{(d{\cdot}i{+}s{-}1)} = \bigcup^{s}_{a{=}0} A_{i}^{(d{\cdot}i{+}a{-}1)},\quad 1 \leq i \leq n,
\end{textstyle}\\[2mm]&\begin{textstyle}
B_{0} = \{\ast\},\quad B_{n} \simeq B.
\end{textstyle}\end{align*}
Then there is a filtration of $E$ by $E\vert_{B_i}$, $0 \leq i \leq n$, as follows (see Whitehead \cite{white2}, for example):
\begin{align*}&\begin{textstyle}
E\vert_{B_{i}} = E\vert_{B_{i{-}1}} \cup_{\Lambda_{i}} C(A_{i}){\times}G,\quad 
\Lambda_i : A_{i}{\times}G \to E\vert_{B_{i{-}1}},\quad 1 \leq i \leq n,
\end{textstyle}\\[2mm]&\begin{textstyle}
E\vert_{B_{0}} = \{\ast\}{\times}G,\quad E\vert_{B_{n}} \simeq E,
\end{textstyle}\end{align*}
and $\tilde{\lambda}_i = \Lambda_i\vert_{A_i} : A_i \to E\vert_{B_{i{-}1}}$ gives a lift of $\lambda_i : A_i \to B_{i-1}$.
Then by induction on $i$, we have 
\begin{align*}&
  E\vert_{B_{i}} = \{\ast\}{\times}G \cup_{\Lambda_{1}} C(A_{1}) {\times} G \cup_{\Lambda_{2}} \cdots \cup_{\Lambda_{i}} C(A_{i}) {\times} G,
\\[2mm]&
\Lambda_{i} : A_{i}{\times}G \xrightarrow{\tilde{\lambda}_{i}{\times}1_{G}}
E\vert_{B_{i-1}}{\times}G
\\[2mm]&\hspace{40pt}= 
\left(
\{\ast\}{\times}G \cup_{\Lambda_{1}} C(A_{1}){\times}G \cdots \cup_{\Lambda_{i{-}1}} C(A_{i{-}1}){\times}G
\right){\times}G
\\[2mm]&\hspace{40pt}
\xrightarrow{1{\times}\mu}
\{\ast\}{\times}G \cup_{\Lambda_{1}} C(A_{1}){\times}G \cdots \cup_{\Lambda_{i{-}1}} C(A_{i{-}1}){\times}G = E\vert_{B_{i-1}},
\end{align*}
where $\mu$ is the multiplication of $G$.
For dimensional reasons, we may regard 
\begin{equation*}
\tilde{\lambda}_{i} : (A_{i},A_{i}^{(d{\cdot}i{+}a{-}1)}) \to (E^{(d{\cdot}i{+}s{-}1)}\vert_{B_{i{-}1}},E^{(d{\cdot}i{+}a{-}1)}\vert_{B_{i{-}1}}),\quad 0 \leq a \leq s,
\end{equation*}
and $\mu(G^{(i)}{\times}G^{(j)}) \subset G^{(i{+}j)}$ up to homotopy.
Then we have the following descriptions for all $k \geq d{\cdot}i{-}1$ and $j \geq d{-}1$:
\begin{align*}&
  E^{(k)}\vert_{B_{i}} = \left(\{\ast\}{\times}G \cup_{\Lambda_{1}} C(A_{1}) {\times} G \cup_{\Lambda_{2}} \cdots \cup_{\Lambda_{i}} C(A_{i}) {\times} G\right)^{(k)},
\\[2mm]&\hspace{32pt}= 
\left(\begin{array}[c]{r}
\{\ast\}{\times}G^{(k)} \cup_{\Lambda_{1}} \bigcup^{s}_{\ell{=}0} (C(A_{1}^{(d{+}\ell{-}1)}) {\times} G^{(k{-}d{-}\ell)}) 
\quad\\\cdots \cup_{\Lambda_{i}} \bigcup^{s}_{\ell{=}0} (C(A_{i}^{(d{\cdot}i{+}\ell{-}1)}) {\times} G^{(k{-}d{\cdot}i{-}\ell)})
\end{array}\right),
\\[2mm]&
\Lambda_{i} : A_{i}^{(d{\cdot}i{+}\ell{-}1)}{\times}G^{(j{-}\ell)} \xrightarrow{\tilde{\lambda}_{i}{\times}1_{G^{(j)}}}
E^{(d{\cdot}i{+}\ell{-}1)}\vert_{B_{i-1}} {\times} G^{(j{-}\ell)}
\\[2mm]&\hspace{20pt}
=\left(\begin{array}[c]{l}
\{\ast\}{\times}G^{(d{\cdot}i{+}\ell{-}1)} 
\\\quad
\cup_{\Lambda_{1}} \bigcup^{s}_{a{=}0} (C(A_{1}^{(d{+}a{-}1)}) {\times} G^{(d{\cdot}(i{-}1){+}\ell{-}a{-}1)})  
\\\quad\cdots 
\cup_{\Lambda_{i{-}1}} \bigcup^{s}_{a{=}0} (C(A_{i{-}1}^{(d{\cdot}(i{-}1){+}a{-}1)}) {\times} G^{(d{+}\ell{-}a{-}1)})
\end{array}\right){\times}G^{(j{-}\ell)}
\\[2mm]&\hspace{18pt}
\xrightarrow{1{\times}\mu}\left(\begin{array}[c]{r}
\{\ast\}{\times}G^{(d{\cdot}i{+}j{-}1)} \cup_{\Lambda_{1}} \bigcup^{s}_{a{=}0} (C(A_{1}^{(d{+}a{-}1)}) {\times} G^{(d{\cdot}(i{-}1){+}j{-}a{-}1)})
\quad\\\cdots
\cup_{\Lambda_{i{-}1}} \bigcup^{s}_{a{=}0} (C(A_{i{-}1}^{(d{\cdot}(i{-}1){+}a{-}1)}) {\times} G^{(d{+}j{-}a{-}1)})
\end{array}\right)
\\[2mm]&\hspace{20pt}= 
\left(\{\ast\}{\times}G \cup_{\Lambda_{1}} C(A_{1}) {\times} G \cup_{\Lambda_{2}} \cdots \cup_{\Lambda_{i{-}1}} C(A_{i{-}1}) {\times} G\right)^{(d{\cdot}i{+}j{-}1)} 
\\[2mm]&\hspace{20pt}= 
E^{(d{\cdot}i{+}j{-}1)}\vert_{B_{i-1}}.
\end{align*}

Similarly, we obtain the following filtration $\{E'_{k}\}_{0 \leq k \leq n{+}m}$ of $E{\times}_GF$.
\begin{align*}&
E'_{k} = 
\left\{\begin{array}[c]{l}
F_{k} \cup_{\Lambda'_{1}} C(A_{1}) {\times} F_{k{-}1} \cup_{\Lambda'_{2}} \cdots \cup_{\Lambda'_{k}} C(A_{k}) {\times} F_{0},\quad k \leq n,
\\[2mm]
F_{k} \cup_{\Lambda'_{1}} C(A_{1}) {\times} F_{k{-}1} \cup_{\Lambda'_{2}} \cdots \cup_{\Lambda'_{n}} C(A_{n}) {\times} F_{k{-}n},\quad n \leq k,
\end{array}\right.
\\[2mm]&
\Lambda'_i : A_{i}{\times}F_{j} 
\xrightarrow{\tilde{\lambda}_i{\times}1_{F_{j}}} E^{(d{\cdot}i{+}s{-}1)}\vert_{B_{i-1}} {\times} F_{j}
\\[2mm]&\hspace{34pt}= 
\left(\begin{array}[c]{r}
G^{(d{\cdot}i{+}s{-}1)} \cup_{\Lambda_{1}} \bigcup^{s}_{a{=}0} (C(A_{1}^{(d{+}a{-}1)}) {\times} G^{(d{\cdot}(i{-}1){+}s{-}a{-}1)}) 
\quad\\\cdots 
\cup_{\Lambda_{i{-}1}} \bigcup^{s}_{a{=}0} (C(A_{i{-}1}^{(d{\cdot}(i{-}1){+}a{-}1)}) {\times} G^{(d{+}s{-}a{-}1)})
\end{array}\right){\times}F_{j} 
\\[2mm]&\hspace{34pt}
\xrightarrow{1{\times}\psi}\left(\begin{array}[c]{r}
F_{i{+}j{-}1} \cup_{\Lambda'_{1}} \bigcup^{s}_{a{=}0} (C(A_{1}^{(d{+}a{-}1)}) {\times} F_{i{+}j{-}2}) 
\quad\\\cdots 
\cup_{\Lambda'_{i{-}1}} \bigcup^{s}_{a{=}0} (C(A_{i{-}1}^{(d{\cdot}(i{-}1){+}a{-}1)}) {\times} F_{j})
\end{array}\right)
\\[2mm]&\hspace{34pt}= 
F_{i{+}j{-}1} \cup_{\Lambda'_{1}} C(A_{1}) {\times} F_{i{+}j{-}2} \cdots \cup_{\Lambda'_{i{-}1}} C(A_{i{-}1}) {\times} F_{j}
\\[2mm]&\hspace{34pt}= 
E'_{i{+}j{-}1}\vert_{B_{i{-}1}},
\end{align*}
since $\psi(G^{(d{\cdot}(\ell{+}1){+}s{-}a{-}1)}{\times}F_{j}) \subseteq \psi(G^{(d{\cdot}(\ell{+}1){+}s{-}1)}{\times}F_{j}) \subset F_{\ell+j}$ by Assumption \ref{ass:compatible-decomposition}.
The above definition of $\Lambda'_{i}$ also determines a map $$\psi_{i,j} : E^{(d{\cdot}(i{+}1){+}s-1)}\vert_{B_i}{\times}F_j \longrightarrow E'_{i+j}\vert_{B_{i}}$$ so that $\Lambda'_{i} = \psi_{i{-}1,j}{\circ}(\tilde{\lambda}_i{\times}1)$.
Let us recall that $F_{j} = F_{j-1} \cup_{\rho_{j}} C(K_{j})$ for $1 \leq j \leq m$.
Then the definition of $E'_{k}$ implies 
\begin{align*}&
    E'_{k} = 
\left\{\begin{array}[c]{l}
\begin{array}[c]{l}
        E'_{k{-}1} \cup C(K_{k}) \cup C(A_{1}){\times}C(K_{k{-}1}) \cup \cdots 
\\
\qquad\quad\cdots \cup C(A_{k{-}1}){\times}C(K_{1}) \cup C(A_{k}){\times}\{\ast\}
\end{array} \hspace{10pt} \textrm{for }k \leq n,
\\[4mm]
\begin{array}[c]{l}
        E'_{k{-}1} \cup C(K_{k}) \cup C(A_{1}){\times}C(K_{k{-}1}) \cup \cdots 
\\
\qquad\quad\cdots \cup C(A_{n{-}1}){\times}C(K_{k{-}n{+}1}) \cup
C(A_{n}){\times}C(K_{k{-}n}) 
\end{array}\hspace{4pt} \textrm{for }k > n.
\end{array}\right.
\end{align*}

To observe the relation between $\Cat{E'_{k{-}1}}$ and $\Cat{E'_{k}}$, we introduce the following two relative homeomorphisms:
\begin{align*}&
\chi(\rho_j) : (C(K_{j}),K_{j}) \to (F_{j{-}1} \cup C(K_{j}),F_{j{-}1}) ~(= (F_{j},F_{j{-}1}))\quad
\\[2mm]&
\chi(\tilde{\lambda}_i) : (C(A_{i}),A_{i}) \to 
\begin{array}[t]{r}
(E^{(d{\cdot}i{+}s{-}1)}\vert_{B_{i}} \cup C(A_{i}),E^{(d{\cdot}i{+}s{-}1)}\vert_{B_{i{-}1}}) \quad\\ (\subset (E^{(d{\cdot}i{+}s)}\vert_{B_{i}},E^{(d{\cdot}i{+}s{-}1)}\vert_{B_{i{-}1}})).
\end{array}
\end{align*}
  Then the attaching map of $C(A_{i}){\times}C(K_{j})$ is given by the Whitehead product $[\chi(\tilde{\lambda}_{i}),\chi(\rho_{j})] : A_{i}{\ast}K_{j} = (C(A_{i}){\times}K_{j}) \cup (A_{i}{\times}C(K_{j})) \to E'_{i{+}j{-}1}$ defined as follows:
\begin{align*}
    & [\chi(\tilde{\lambda}_{i}),\chi(\rho_{j})]\vert_{C(A_{i}){\times}K_{j}}
    :
        \begin{array}[t]{l}
        C(A_{i}){\times}K_{j} 
        \xrightarrow{\chi(\tilde{\lambda}_{i}){\times}1} 
        E^{(d{\cdot}i{+}s)}\vert_{B_{i}} {\times} F_{j{-}1} 
        \\
        \hspace{-20pt}\subseteq E^{(d{\cdot}(i{+}1){+}s{-}1)}\vert_{B_{i}} {\times} F_{j{-}1} 
        \xrightarrow{\psi_{i,j{-}1}}
    E'_{i{+}j{-}1}\vert_{B_{i}} \subseteq E'_{i{+}j{-}1},
        \end{array}
\\[2mm]
    & [\chi(\tilde{\lambda}_{i}),\chi(\rho_{j})]\vert_{A_{i}{\times}C(K_{j})} : A_{i}{\times}C(K_{j}) 
        \xrightarrow{\tilde{\lambda}_{i}{\times}\chi(\rho_{j})}
        \begin{array}[t]{l}
        E^{(d{\cdot}i{+}s{-}1)}\vert_{B_{i{-}1}} {\times} F_{j} 
        \\
        \hspace{-10pt}\xrightarrow{\psi_{i{-}1,j}} E'_{i{+}j{-}1}\vert_{B_{i{-}1}} \subseteq E'_{i{+}j{-}1}.
        \end{array}
\end{align*}
This implies immediately that $\Cat{E'_{k}} \leq \Cat{E'_{k{-}1}}+1$.
Then by induction on $k$, we obtain that 
$\Cat{E'_{k}} \leq k$.
  Thus we have $\begin{textstyle}\Cat{X} = \Cat{E{\times}_GF} = \Cat{E'_{m{+}n}} \leq m{+}n \leq m{+}\frac{\dim B}{d}\end{textstyle}$.
  This completes the proof of Theorem \ref{thm:main}.

\section{Proof of Theorem \ref{thm:main1}}\label{sect:3}

As is well known, we have the following principal bundles (see for
example \cite{adams}, \cite{yok4} and \cite{harvey} in particular for the last fibration):
\begin{align*}
  & \mathrm{Sp}(1) \longrightarrow \mathrm{Sp}(2) \longrightarrow S^7, \\
  & \mathrm{SU}(3) \longrightarrow \mathrm{SU}(4) \longrightarrow S^7, \\
  & G_2 \longrightarrow \mathrm{Spin}(7) \longrightarrow S^7, \\
  & \mathrm{Spin}(7) \longrightarrow \mathrm{Spin}(9) \longrightarrow S^{15}, \\
  & G_2 \longrightarrow \mathrm{Spin}(8) \longrightarrow S^7 \times S^7.
\end{align*}
Each scalar matrix $(-1) \in \mathrm{Sp}(2)$ and $(-1) \in \mathrm{SU}(4)$ acts on $S^7$
as the antipodal map, and so does the center of $\mathrm{Spin}(7)$.
Similarly the center of $\mathrm{Spin}(9)$ acts on $S^{15}$ as the antipodal map.
Recall that the center of $\mathrm{Spin}(8)$ is isomorphic to $\mathbb Z/2 \times
\mathbb Z/2$, each generator of which acts on $S^7$ as the antipodal map
respectively.
Since there are isomorphisms $\mathrm{Sp}(2) \cong
\mathrm{Spin}(5)$ and $\mathrm{SU}(4) \cong \mathrm{Spin}(6)$, we obtain principal bundles:
\begin{align*}
  & \mathrm{Sp}(1) \longrightarrow \mathrm{SO}(5) \longrightarrow \mathbb RP^7, \\
  & \mathrm{SU}(3) \longrightarrow \mathrm{SO}(6) \longrightarrow \mathbb RP^7, \\
  & \mathrm{G}_2 \longrightarrow \mathrm{SO}(7) \longrightarrow \mathbb RP^7, \\
  & \mathrm{Spin}(7) \longrightarrow \mathrm{SO}(9) \longrightarrow \mathbb RP^{15}, \\
  & \mathrm{G}_2 \longrightarrow \mathrm{PO}(8) \longrightarrow \mathbb RP^7 \times
  \mathbb RP^7.
\end{align*}

\noindent
Cone decompositions of the fibres except $\mathrm{Spin}(7)$ are given as follows (see Theorem 2.1 of \cite{IM:cat-sp(3)} for $\mathrm{G}_{2}$):
  \begin{align*}
    & \ast \subset \mathrm{Sp}(1) = S^3, \\
    & \ast \subset \mathrm{SU}(3)^{(5)} \subset \mathrm{SU}(3), \\
    & \ast \subset \mathrm{G}_2^{(5)} \subset \mathrm{G}_2^{(8)} \subset \mathrm{G}_2^{(11)} \subset \mathrm{G}_2,
  \end{align*}
  where $\mathrm{SU}(3)^{(5)} = \mathrm{G}_2^{(5)} = \Sigma \mathbb CP^2$,
  $\mathrm{SU}(3) = \mathrm{SU}(3)^{(5)} \cup CS^7$,
  $\mathrm{G}_2^{(8)} \simeq \mathrm{G}_2^{(5)} \cup C(S^5 \cup e^7)$,
  $\mathrm{G}_2^{(11)} \simeq \mathrm{G}_2^{(8)} \cup C(S^8 \cup e^{10})$ and
  $\mathrm{G}_2 = \mathrm{G}_2^{(11)} \cup CS^{13}$.
  Since these fibres satisfy the conditions in Remark
  \ref{rem:cellular} (\ref{rem:cellular1}), we obtain $\Cat{\mathrm{SO}(5)} \leq 8$, $\Cat{\mathrm{SO}(6)} \leq 9$,
  $\Cat{\mathrm{SO}(7)} \leq 11$ and $\Cat{\mathrm{PO}(8)} \leq 18$ using Theorem \ref{thm:main}.
  By virtue of the mod 2 cup-lengths we have that $\cuplen(\mathrm{SO}(5)) \geq
  8$, $\cuplen(\mathrm{SO}(6)) \geq 9$, $\cuplen(\mathrm{SO}(7)) \geq 11$ and
  $\cuplen(\mathrm{PO}(8)) \geq 18$ respectively.
  Thus we obtain the results for $\mathrm{SO}(5)$, $\mathrm{SO}(6)$, $\mathrm{SO}(7)$ and $\mathrm{PO}(8)$.

A cone decomposition of $\mathrm{Spin}(7)$ is given as follows in \cite{IMN:spin(7)}:
  \begin{align*}
    \ast = F_0 \subset F_1 \subset F_2 \subset F_3 \subset F_4 \subset
    F_5 = \mathrm{Spin}(7),
  \end{align*}
  where $F_1 = \mathrm{SU}(4)^{(7)}$, $F_2 = \mathrm{SU}(4)^{(12)} \cup e^6$,
  $F_3 = \mathrm{SU}(4) \cup e^6 \cup e^9 \cup e^{11} \cup e^{13}$ and
  $F_4 = \mathrm{Spin}(7)^{(18)}$.
  We need here to check if the filtration satisfies Assumption 1; the only problem is to determine whether $\psi \vert_{\mathrm{Spin}(7)^{(3)}{\times}F_1} : \mathrm{Spin}(7)^{(3)}{\times}F_{1} \to F$ is compressible into $F_{4}$ or not.
  Since $\mathrm{Spin}(7)^{(3)}$ and $F_1$ are included in $\mathrm{SU}(4) \subset F_4$, we have ${\rm Im}\, (\psi \vert_{\mathrm{Spin}(7)^{(3)}{\times}F_1}) \subset F_4$.
  Then we obtain $\Cat{\mathrm{SO}(9)} \leq 20$ using Theorem \ref{thm:main}.
  The mod 2 cup-length implies that $\cuplen(\mathrm{SO}(9)) \geq 20$.
  Thus we obtain the result for $\mathrm{SO}(9)$.

Since $\mathrm{SO}(8)$ is homeomorphic to $\mathrm{SO}(7) \times S^7$,
we easily see that
\begin{displaymath}
  \Cat{\mathrm{SO}(8)} \leq \Cat{\mathrm{SO}(7)}+\Cat{S^7}=12
\end{displaymath}
by Takens \cite{takens}.
  The mod 2 cup-length implies that $\cuplen(\mathrm{SO}(8)) \geq 12$.
  Thus we obtain the result for $\mathrm{SO}(8)$.
  This completes the proof of Theorem \ref{thm:main1}.

\section{Proof of Theorems \ref{thm:maing} and \ref{thm:main2}}\label{sect:4}

Firstly, we show Theorem \ref{thm:maing}.
The following principal bundle is well-known:
\begin{equation*}
  \mathrm{SU}(n{-}1) \longrightarrow \mathrm{SU}(n) \longrightarrow S^{2n-1}.
\end{equation*}
The central (cyclic) subgroup $C_m$ of $\mathrm{SU}(n)$ acts on $S^{2n-1}$
freely and hence we obtain a principal bundle:
\begin{equation*}
  \mathrm{SU}(n{-}1) \longrightarrow \mathrm{SU}(n)/C_m
  \longrightarrow L^{2n-1}(m),
\end{equation*}
where $L^{2n-1}(m)$ is a lens space of dimension $2n{-}1$.

A cone decomposition of $\mathrm{SU}(n{-}1)$ is constructed by
Kadzisa \cite{Kadzisa:cone-decomp-su}:
  \begin{displaymath}
    \ast \subset V \subset V^2 \subset \cdots \subset V^{n{-}2} = SU(n{-}1),
  \end{displaymath}
where $V^k \subseteq SU(n{-}1)$ is a representing subspace of the quotient module $H^*(SU(n{-}1))/D^{k+1}$ and $D^{k+1}$ is the submodule generated by products of $k{+}1$ elements in positive degrees, which satisfies $V^i{\cdot}V^j \subseteq V^{i+j}$ for any $i$ and $j$.
Thus $V$ is the subcomplex $S^3 \cup e^5 \cup e^7 \cup \cdots \cup e^{2n{-}3}$ of $\mathrm{SU}(n{-}1)$ which is homeomorphic to $\Sigma\mathbb CP^{n{-}2}$ (see \cite{yok1}, for example).
Then Assumption 1 is automatically satisfied, and hence using $\mathrm{SU}(n{-}1)^{(k)} \subset V^k$, we obtain 
  \begin{displaymath}
    \Cat{\mathrm{SU}(n)/C_m} \leq 3(n{-}1)
  \end{displaymath}
by Theorem \ref{thm:main}.
This completes the proof of Theorem \ref{thm:maing}.

Secondly, we show Theorem \ref{thm:main2}.
By Rudyak \cite{Rudyak:ls-cat_mfds1} \cite{Rudyak:ls-cat_mfds2} and Strom \cite{strom}, we know the following proposition.
\begin{prop}
  \label{prop:weight}
  {\rm (Rudyak \cite{Rudyak:ls-cat_mfds1} \cite{Rudyak:ls-cat_mfds2}, Strom \cite{strom})}
Let $h$ be a cohomology theory.
For an element $u \in \tilde{h}^*(X)$, let $\weight(u;h)$ be the
minimal number $k$ such that ${(e^X_k)}^*(u) \not= 0$ where $e^X_k : P^k\Omega{X} \to P^{\infty}\Omega X \simeq X$, which satisfies
  \begin{enumerate}
  \item\label{wgt(1)} We have $\weight(0;h) = \infty$ and $\infty > \weight(u;h) \geq 1$ for any $u \not= 0$ in $\tilde{h}^*(X)$.
  \item\label{wgt(2)} For any cohomology theory $h$, we have
  \begin{equation*}
    \min\left\{\weight(u;h),\weight(v;h)\right\} \leq \weight(u+v;h).
  \end{equation*}
  \item\label{wgt(3)} For any multiplicative cohomology theory $h$, we have
  \begin{equation*}
    \weight(u;h) + \weight(v;h) \leq \weight(u{\cdot}v;h).
  \end{equation*}
  \item\label{wgt(4)} $\weight(X;h) = \max \{ \weight(u;h) \mid u \in \tilde{h}^*(X), u\not= 0\}$.
  \end{enumerate}
\end{prop}

Le us recall that, for any compact Lie group $G$, the ordinary cohomology of $\Omega G$ is concentrated in even degrees.
Then, for any element $u$ of even degree in $\tilde{H}^*(G;\mathbb Z/p)$, we have $\weight(u;{H\mathbb Z/p}) \geq 2$, since $P^1(\Omega G) = \Sigma \Omega (G)$.

The cohomology rings of $\mathrm{SU}(p^r)/C_{p^s}$ for a prime $p$ and
$1 \leq s \leq r$ are given as follows (see \cite{baumbr}):
\begin{displaymath}
  H^*(\mathrm{SU}(p^r)/C_{p^s};\mathbb Z/p) = \mathbb Z/p[x_2]/(x_2^{p^r})
  \otimes {\wedge} (x_1, x_3, \ldots, x_{2p^r-3}).
\end{displaymath}
Note that $x_1^2 = x_2$ if $p=2$ and $s=1$.
Then, using Proposition \ref{prop:weight}, we obtain
\begin{displaymath}
  \weight(\mathrm{SU}(p^r)/C_{p^s};{H\mathbb Z/p}) \geq \weight
  (x_1{\cdot}x_2^{p^r-1}{\cdot}x_3 \cdot \cdots \cdot x_{2p^r-3};{H\mathbb Z/p})
  \geq 3(p^r-1),
\end{displaymath}
since $\weight(x_2;{H\mathbb Z/p}) \geq 2$.
Thus we have the following lemma.
\begin{lem}
  $\sigmacat{\mathrm{SU}(p^r)/C_{p^s}} \geq 3(p^r{-}1)$ for any prime
  $p$ and $1 \leq s \leq r$.
\end{lem}
By using Theorem \ref{thm:maing}, we obtain Theorem \ref{thm:main2}.
%
%

\end{document}